\documentclass[a4paper,12pt,twoside]{article}
\usepackage[english]{babel}
\usepackage[T1]{fontenc}
\usepackage{amsmath,amssymb,amsfonts,amsthm,amsopn,amstext,amscd,latexsym,enumerate}
\usepackage{graphicx}
\usepackage{float}
\usepackage{makeidx}
\usepackage{algorithm}
\usepackage{algorithmic}
\usepackage{slashbox}

\catcode `\Á=13 \catcode `\É=13 \catcode `\Í=13 \catcode `\Ó=13
\catcode `\Ú=13

\title{{ Extensions of linear regression models based on set arithmetic
for interval data}}

\author{\hspace*{-2cm} Angela Blanco-Fern\'{a}ndez$^{\rm a}$\thanks{$^\ast$Corresponding author. Email: blancoangela@uniovi.es
}, Marta Garc\'{i}a-B\'{a}rzana$^{\rm a}$, Ana Colubi$^{\rm a}$ and\\ \hspace*{-5cm}Erricos J. Kontoghiorghes$^{\rm b}$ \\  \hspace*{-2cm}$^{\rm a}${\em{\small{Department of Statistics, University of Oviedo, Calvo Sotelo s/n, Oviedo 33007, Spain}}};\\\hspace*{-3cm} $^{\rm b}${\em{\small{Department of Commerce, Finance and Shipping, Cyprus University of Technology}}}\\ {\em{\small{140 Ayiou Andreou Street, Lemmesos, Cyprus}}} }
\date{}
\begin{document}

\maketitle

\begin{abstract}
Extensions of previous linear regression models for interval data  are
presented. A more flexible
simple linear model is formalized. The new model may express
cross-relationships between mid-points
and spreads of the interval data in a unique equation based on the interval arithmetic. Moreover, extensions
to the multiple case are addressed. The associated least-squares estimation problems
are solved. Empirical results and a real-life
application are
presented in order to show the applicability and the differences among the
proposed models.

{\bf keywords:} multiple linear regression model; interval data; set arithmetic;
least-squares estimation

\end{abstract}

\section{Introduction}
The statistical treatment of interval data is recently being considered extensively (see \cite{BeresteanuMolinari:08, Yuetal:06,DursoGiordani:06, Giordani:10, Hashimotoetal:11, Tseetal:08}).
Interval data are useful to model variables with uncertainty in their formalization, due to an
imprecise observation or an inexact measurement,
fluctuations, grouped data or censoring. Linear regression models for interval data have been
previously analyzed (see \cite{Diamond:90,Giletal:02, GRetal:07, GRetal:07b, Blancoetal:11, Blancoetal:12,BillardyDiday:00, LimaNetoetal:10, ManskiTamer:02}). Regression models with interval-valued explanatory variables and
interval-valued response are considered. There are two main approaches to face these kinds of
problems. One is based on fitting separate models for mid-points and spreads
(see \cite{BillardyDiday:00,LimaNetoetal:10}). This approach has not been considered  under
probabilistic assumptions on the population models, and inferential studies have not been
developed yet. This is non-trivial, since the non-negativity constraints satisfied by the
spread variables prevent the corresponding model to be treated as a classical linear regression.
Thus, although the usual
fitting techniques are used, the associated inferences are no longer valid. The second
approach overcomes this difficulty by considering a model based on the set arithmetic
(see \cite{GRetal:07,Blancoetal:11}). The least squares estimators are found as
solutions of constrained minimization problems and inferential studies have been developed in \cite{GRetal:07b} and \cite{Blancoetal:12}, among others.  

Extensions for the simple linear regression models within the framework of the work in \cite{GRetal:07} and \cite{Blancoetal:11} are developed. On one hand, a more flexible simple linear model is formalized. The previous regression functions model the response mid-points (respectively spreads)
by means of the explanatory mid-points (respectively spreads). The new
model is able to accommodate cross-relationships between mid-points
and spreads in a unique equation based on the set arithmetic. As the model in \cite{Blancoetal:11},
the new one is based on the so-called \textit{canonical decomposition} of the intervals. 
On the other hand, extensions
to the multiple case are addressed. Due to the essential differences of the model
in \cite{GRetal:07} and those based on the canonical decomposition, two multiple models will be
introduced.
The least-squares (LS) estimation problems associated with the
proposed regression models are solved. Some empirical results and a real-life
application are presented in order to show the applicability and the differences among the
proposed models.

The rest of the paper is organized as follows: In Section
\ref{section-preliminaries} some preliminary concepts about the interval
framework are presented and several previous simple
linear models based on the set arithmetic are
revised. Extensions of those linear models are introduced in Sections
\ref{subsection-M-MG}, \ref{subsection-basic-multiple} and \ref{subsection-MG-multiple}. 
The theoretical formalization and the associated LS 
estimation problems
are addressed. In Section \ref{section-empirical} the empirical
performance and the practical applicability of the models are shown through some
simulation studies and a real-life case-study. Finally, Section
\ref{section-conclusions} includes some conclusions and future directions.

\section{Preliminaries}\label{section-preliminaries}

The considered interval experimental data are elements
belonging to the space
$\mathcal{K _{\mathrm{c}}}(\mathbb{R})=\{[a_1,a_2] : a_1,a_2 \in \mathbb{R}, a_1
\le a_2\}$.
Each interval $A\in \mathcal{K _{\mathrm{c}}}(\mathbb{R}) $ can be parametrized in terms of its
mid-point,
$\textrm{mid}\ \hspace*{-0.1cm}  A =(\sup A + \inf A)/2$, and its spread, $\textrm{spr}\
\hspace*{-0.1cm}  A = (\sup A - \inf A)/2$. The notation
$A=[\mathrm{mid}A\pm \mathrm{spr} A]$ will be used.  An
alternative representation for intervals is the so-called {\it canonical
decomposition}, introduced in \cite{Blancoetal:11}, given by $A=\mathrm{mid}A[1\pm 0]
+ \mathrm{spr}A[0\pm 1]$. It allows the consideration of the {\it mid} and
{\it spr} components of $A$ separately within the interval arithmetic.

The Minkowski addition and the product by scalars form the natural arithmetic on $\mathcal{K
_{\mathrm{c}}}(\mathbb{R})$. In terms of the (mid, spr)-representation these operations
can be
jointly expressed as $$A+ \lambda B=[( \mathrm{mid}A + \lambda\mathrm{mid}B ) \
\pm \ (\mathrm{spr} A + |\lambda|\mathrm{spr} B)]$$ for                      
any $A,B \in \mathcal{K _{\mathrm{c}}}(\mathbb{R})$ and $\lambda \in
\mathbb{R}$. The space $(\mathcal{K
_{\mathrm{c}}}(\mathbb{R}),+,\ldotp)$ is not linear but semilinear (or conical), due to the
lack of symmetric element with respect to the addition.  $ \mathcal{K
_{\mathrm{c}}}(\mathbb{R})$ can be identified with the cone $\mathbb{R}\times\mathbb{R}^+$ of
$\mathbb{R}^2$.  The expression $A+(-1)B$
generally differs from the natural difference $A-B$. 
If it exists $C=A-_H B \in \mathcal{K
_{\mathrm{c}}}(\mathbb{R})$ verifying that $A=B+C$, $C$ is called {\it Hukuhara difference}
between the pair of intervals $A$ and $B$. 
The interval $C$ exists iff $\mathrm{spr} B \le \mathrm{spr} A$. 

For every $A,B\in \mathcal{K _{\mathrm{c}}}(\mathbb{R})$, the $L_2$-type generic metric in
\cite{Trutschnigetal:09} is defined as
\begin{equation}\label{dtita}
d_{\theta}(A,B)=((\mathrm{mid}A-\mathrm{mid}B)^2+\theta\,(\mathrm{spr}A-\mathrm{
spr}B)^2)^{\frac{1}{2}}
\end{equation}
\noindent for an arbitrary $\theta \in (0,\infty)$. 

Given a probability space
 $(\Omega, \mathcal {A}, P)$, the mapping $\emph{\textbf{x}}:\Omega \rightarrow
\mathcal{K _{\mathrm{c}}}(\mathbb{R})$ is said to be a {\it random interval}
iff
$\mathrm{mid}\,\emph{\textbf{x}}, \mathrm{spr}\,\emph{\textbf{x}} : \Omega
\rightarrow \mathbb{R}$ are real random variables and
$\mathrm{spr}\,\emph{\textbf{x}} \ge 0$.
Random intervals will be denoted with bold lowercase letters,
$\emph{\textbf{x}}$, random interval-valued vectors will be represented by
non-bold lowercase letters, $x$, and interval-valued matrices will be denoted
with uppercase letters, $X$.

The expected value of
${\emph{\textbf{x}}}$ is defined in terms of the well-known Aumann expectation,
which satisfies that
\begin{equation}\label{aumann}
E(\emph{\textbf{x}})=[E(\mathrm{mid}\,\emph{\textbf{x}}) \pm
E(\mathrm{spr}\,\emph{\textbf{x}})] ,
\end{equation}
\noindent whenever 
$\mathrm{mid}\,\emph{\textbf{x}}, \mathrm{spr}\,\emph{\textbf{x}} \in L^1$. The variance of a random
interval
$\emph{\textbf{x}}$ can be defined as the usual {\it Fr\'{e}chet variance} (see \cite{Nather:97}) 
associated with the Aumann expectation in the metric space $(\mathcal{K
_{\mathrm{c}}}(\mathbb{R}), d_\theta)$, i.e.
\begin{equation}\label{variance}
\nonumber\sigma^2_{\emph{\textbf{x}}}=E\Big(d_\theta^2(\emph{\textbf{x}},E(\emph
{\textbf{x}}))\Big)\,=\sigma^2_{\textrm{mid}
\emph{\textbf{x}}}+\theta\sigma^2_{\textrm{spr} \emph{\textbf{x}}},
\end{equation}
\noindent whenever $\mathrm{mid}\,\emph{\textbf{x}}, \mathrm{spr}\,\emph{\textbf{x}} \in L^2$. The conical structure
of the space $\mathcal{K _{\mathrm{c}}}(\mathbb{R})$ entails some differences
to define the usual covariance (see \cite{Korner:97}). In terms of the $d_\theta$
metric it has the
expression
\begin{eqnarray}\label{cov}
\nonumber
\sigma_{\emph{\textbf{x}},\emph{\textbf{y}}}=\sigma_{\textrm{mid}
\emph{\textbf{x}},\textrm{mid} \emph{\textbf{y}} }+ \theta\sigma_{\textrm{spr}
\emph{\textbf{x}},\textrm{spr} \emph{\textbf{y}}},
\end{eqnarray}
 whenever those classical covariances exist. The expression
$Cov(x,y)$ denotes the
covariance matrix between two random interval-valued vectors
$x=(\emph{\textbf{x}}_1,\ldots,\emph{\textbf{x}}_k)$ and 
$y=(\emph{\textbf{y}}_1,\ldots,\emph{\textbf{y}}_k)$.

Let $\emph{\textbf{x}},\emph{\textbf{y}} : \Omega\rightarrow \mathcal{K
_{\mathrm{c}}}(\mathbb{R})$ be two random intervals. The {\it basic} simple linear model (see \cite{Giletal:02}) to relate two random
intervals has the expression:
\begin{equation}\label{model-basic}
\emph{\textbf{y}}=b \emph{\textbf{x}}+\boldsymbol{\varepsilon}
\end{equation}
with $b\in \mathbb{R}$ and $\boldsymbol{\varepsilon}:\Omega \rightarrow \mathcal{K
_{\mathrm{c}}}(\mathbb{R})$ is an interval-valued random error such that
$E[\boldsymbol{\varepsilon} | \emph{\textbf{x}}]=\Delta \in \mathcal{K
_{\mathrm{c}}}(\mathbb{R})$. The  LS
estimation of (\ref{model-basic}) has been solved analytically by
means of a constrained minimization problem in \cite{GRetal:07}.

Model (\ref{model-basic}) only involves one regression parameter $b$ to model the dependency.
Thus, it induces quite restrictive separate models for the {\it mid} and {\it spr} components of the intervals. Specifically,  $\mathrm{mid}
\emph{\textbf{y}}=b
\mathrm{mid} \emph{\textbf{x}} + \mathrm{mid}\boldsymbol{\varepsilon}$  and
$\mathrm{spr} \emph{\textbf{y}}=|b| \mathrm{spr} \emph{\textbf{x}} +
\mathrm{spr}\boldsymbol{\varepsilon}$. 

A more flexible linear model, called model M, has been
introduced in \cite{Blancoetal:11}. It is
defined in terms of the {\it canonical decomposition} as follows:
\begin{equation}\label{model-M}
\emph{\textbf{y}}=b^1 \mathrm{mid}\emph{\textbf{x}}\,[1 \pm
0]+b^2 \mathrm{spr}\emph{\textbf{x}}\,[0 \pm 1]+\gamma\,[1 \pm
0]+\boldsymbol{\varepsilon},
\end{equation}
\noindent where $b^1, b^2\in \mathbb{R}$ are the regression coefficients,
$\gamma\in \mathbb{R}$ is an intercept term influencing the {\it mid} component
of $\emph{\textbf{y}}$ and $\boldsymbol{\varepsilon}$ is a random interval error
satisfying that $E[\boldsymbol\varepsilon|\emph{\textbf{x}}]=[-\delta,\delta]$, with
$\delta \ge 0$. From
(\ref{model-M}) the linear relationships $\mathrm{mid} \emph{\textbf{y}}=b^1
\mathrm{mid} \emph{\textbf{x}} + \gamma + \mathrm{mid}\boldsymbol{\varepsilon}$  and
$\mathrm{spr} \emph{\textbf{y}}=|b^2| \mathrm{spr} \emph{\textbf{x}} +
\mathrm{spr}\boldsymbol{\varepsilon}$ are transferred, where $b^1$ and $b^2$ may be different.  
The LS  estimation  leads to analytic
expressions of the regression parameters of model M (see \cite{Blancoetal:11}). Confidence sets based
on those estimators have been developed in \cite{Blancoetal:12}.

\section{A flexible simple linear regression model: the model 
$\mathrm{M}_G$}\label{subsection-M-MG}

Following (\ref{model-M}), the
 model $\mathrm{M}_G$ between $\emph{\textbf{x}}$ and
$\emph{\textbf{y}}$ is defined as: 
\begin{equation}\label{model-Mg}
\noindent  \emph{\textbf{y}}\hspace*{-0.05cm}=\hspace*{-0.05cm}b^1\mathrm{mid} \emph{\textbf{x}}
[1\pm 0] + b^2\mathrm{spr} \emph{\textbf{x}} [0 \pm 1] +b^3\mathrm{mid}
\emph{\textbf{x}} [0\pm 1] +b^4\mathrm{spr} \emph{\textbf{x}} [1\pm 0] + \gamma
[1\pm 0] + \boldsymbol\varepsilon,
\end{equation}
where $b^i, \gamma \in \mathbb{R}$, $i=1,\ldots, 4$ and  $E(\boldsymbol\varepsilon \vert
 \emph{\textbf{x}})=[-\delta,\delta]\in \mathcal{K}_c(\mathbb{R})$, $\delta \ge 0$.
The linear relationships for the {\it mid} and {\it spr} variables
transferred from (\ref{model-Mg}) are
$$\mathrm{mid} \emph{\textbf{y}}= b^1\mathrm{mid} \emph{\textbf{x}} +
b^4\mathrm{spr} \emph{\textbf{x}} + \gamma + \mathrm{mid} \boldsymbol\varepsilon$$
\vspace*{-0.9cm}

\noindent and
\vspace*{-0.7cm}

$$\mathrm{spr} \emph{\textbf{y}}=  |b^2|\mathrm{spr} \emph{\textbf{x}} +
|b^3||\mathrm{mid} \emph{\textbf{x}}| + \mathrm{spr} \boldsymbol\varepsilon .$$
Thus, both variables $\mathrm{mid} \emph{\textbf{y}}$ and
$\mathrm{spr} \emph{\textbf{y}}$ are modelled from the complete information
provided by the independent random interval $\emph{\textbf{x}}$, characterized by the
random vector $(\mathrm{mid} \emph{\textbf{x}},\mathrm{spr} \emph{\textbf{x}})$.

For a simpler notation, the random intervals defined from
$\emph{\textbf{x}}$ are denoted by $\emph{\textbf{x}}^M$, $\emph{\textbf{x}}^S$,
$\emph{\textbf{x}}^C$ and $\emph{\textbf{x}}^R$, in the same order as they
appear in (\ref{model-Mg}). Thus, the model $\mathrm{M}_G$ is equivalently
expressed as: 
\begin{equation}\label{model-Mg2}
\nonumber\emph{\textbf{y}}=b^1\ \hspace*{-0.1cm} \emph{\textbf{x}}^M+b^2 \
\hspace*{-0.1cm} \emph{\textbf{x}}^S+b^3 \ \hspace*{-0.1cm}
\emph{\textbf{x}}^C+b^4\ \hspace*{-0.1cm} \emph{\textbf{x}}^R + \gamma[1\pm 0] +
\boldsymbol\varepsilon.
\end{equation}

\noindent Moreover, in order to unify the notation for the estimation problems
of the different linear models,  the real interval
$\Delta=[\gamma-\delta,\gamma+\delta]$ is defined. Then, the regression function
associated with the model $\mathrm{M}_G$ can be written as: 
\begin{equation}\label{regression_function_MG}
E(\emph{\textbf{y}} | \emph{\textbf{x}})=b^1\ \hspace*{-0.1cm}
\emph{\textbf{x}}^M+b^2 \ \hspace*{-0.1cm} \emph{\textbf{x}}^S+b^3 \
\hspace*{-0.1cm} \emph{\textbf{x}}^C+b^4\ \hspace*{-0.1cm} \emph{\textbf{x}}^R +
\Delta .
\end{equation}

Since
$\emph{\textbf{x}}^S=-\emph{\textbf{x}}^S$ and
$\emph{\textbf{x}}^C=-\emph{\textbf{x}}^C$, the model $\mathrm{M}_G$ always admits four equivalent
expressions. This
property allows the simplification of the estimation process, because it is possible to
search only for non-negative estimates of the parameters $b^2$
and $b^3$. 

Given a random sample
$\{\left(\emph{\textbf{x}}_{j},\emph{\textbf{y}}_j\right)\}_{j=1}^{n}$ obtained
from two random intervals $(\emph{\textbf{x}},\emph{\textbf{y}})$ verifying 
(\ref{model-Mg}), the LS estimation of the parameters $(b^1,b^2,b^3,b^4,\Delta)$
in  (\ref{regression_function_MG}) consists in minimizing 
\begin{equation}\label{minim-Mg}
 \frac{1}{n} \sum_{i=1}^{n}
d_\theta^{2}(\emph{\textbf{y}}_i,a\emph{\textbf{x}}_i^M+b\emph{\textbf{x}}
_i^S+c\emph{\textbf{x}}_i^C+d\emph{\textbf{x}}_i^R+C) 
\end{equation}
over
$(a,b,c,d,C)\in\mathbb{R}^4\times \mathcal{K}_c(\mathbb{R})$.
However, since from equation (\ref{model-Mg}), $\boldsymbol\varepsilon_i=\emph{\textbf{y}}_i -_H
(b^1\emph{\textbf{x}}_i^M+b^2\emph{\textbf{x}}_i^S+b^3\emph{\textbf{x}}_i^C+b^4\emph{
\textbf{x}}_i^R)$, (\ref{minim-Mg}) must be
optimized over a suitable feasible set assuring the existence of the sample
residuals, i.e., the corresponding Hukuhara differences.
Note that
$$\mathrm{spr}(a\emph{\textbf{x}}_i^M+b\emph{\textbf{x}}_i^S+c\emph{\textbf{x}}
_i^C+d\emph{\textbf{x}}_i^R)=|b| \mathrm{spr} \emph{\textbf{x}}_i +|c|
|\textrm{mid} \emph{\textbf{x}}_i|$$ for all $i=1,\ldots,n$ and 
$b^2$ and $b^3$ can be assumed to be non-negative. Then, taking into account the condition guaranteeing
the existence of the Hukuhara difference, the feasible set can be expressed as
\begin{equation}\label{SG}
\hspace*{-0.2cm} \Gamma_{G}\hspace*{-0.1cm}=\hspace*{-0.1cm}\{(b,c)\in
[0,\infty)\hspace*{-0.1cm} \times \hspace*{-0.1cm}[0,\infty)\hspace*{-0.1cm}:
\hspace*{-0.1cm}b\ \hspace*{-0.05cm}\mathrm{spr} \emph{\textbf{x}}_i + c \
\hspace*{-0.1cm}|\textrm{mid} \emph{\textbf{x}}_i| \le \mathrm{spr}
\emph{\textbf{y}}_i , \forall i=1,\ldots,n\}.
\end{equation}

 If 
$(\widehat{b^1},\widehat{b^2},\widehat{b^3},\widehat{b^4})\in \mathbb{R}^4$ denotes a feasible
estimate, then
the interval parameter $\Delta$ can be directly estimated by 
$$\widehat{\Delta}=\overline{\emph{\textbf{y}}} -_H  \left(\widehat{b^1}
\overline{\emph{\textbf{x}}^M}+\widehat{b^2}\overline{\emph{\textbf{x}}^S}
+\widehat{b^3}\overline{\emph{\textbf{x}}^C}+\widehat{b^4}\overline{\emph{
\textbf{x}}^R}\right).$$

\noindent As a result, the LS minimization problem is 
\begin{equation}\label{minim-Mg-2} 
\hspace*{-0.15cm}\min_{{\small \begin{array}{c}(a,d)\in \mathbb{R}^2\\(b,c)\in
\Gamma_G\end{array}}}\hspace*{-0.15cm}\frac{1}{n} \hspace*{-0.1cm}\sum_{i=1}^{n}
\hspace*{-0.1cm}d_\theta^{2}\hspace*{-0.05cm}\Big(\emph{\textbf{y}}_i-_H
(a\emph{\textbf{x}}_i^M+b\emph{\textbf{x}}_i^S+c\emph{\textbf{x}}_i^C+d\emph{
\textbf{x}}_i^R), \overline{\emph{\textbf{y}}} -_H (a
\overline{\emph{\textbf{x}}^M}+b\overline{\emph{\textbf{x}}^S}+c\overline{\emph{
\textbf{x}}^C}+d\overline{\emph{\textbf{x}}^R})\hspace*{-0.1cm}\Big).
\end{equation}

\noindent The problem (\ref{minim-Mg-2}) can be solved separately for $(a,d)$ and $(b,c)$. The
minimization over $(a,d)$ is done without restrictions and it leads to
the following analytic estimators of $(b^1,b^4)$ in the model 
$\mathrm{M}_G$:
\begin{equation}\label{estimator-MG-b1b4}
(\widehat{b^1},\widehat{b^4})^t=S_1^{-1} z_1.
\end{equation}
\noindent Here $z_1=(\widehat{\sigma}_{\emph{\textbf{x}}^M,\emph{\textbf{y}}}  
 ,  \widehat{\sigma}_{\emph{\textbf{x}}^R,\emph{\textbf{y}}})^t$ and $S_1$
corresponds to the sample covariance matrix of the interval-valued random vector
$(\emph{\textbf{x}}^M,\emph{\textbf{x}}^R)$.

The minimization for $(b,c)$ is performed over the
feasible set  $\Gamma_G$, which is nonempty, closed and convex. The objective function to
be minimized over $(b,c)$ can be
expressed as the globally convex function
\begin{equation}\label{minim-MG-b2b3}
g(b,c)=(b)^2\ \hspace*{-0.1cm}\widehat{\sigma}^2_{\emph{\textbf{x}}^S} + (c)^2\
\hspace*{-0.1cm}\widehat{\sigma}^2_{\emph{\textbf{x}}^C}  +2\ \hspace*{-0.1cm}b\
\hspace*{-0.1cm}c\
\hspace*{-0.1cm}\widehat{\sigma}_{\emph{\textbf{x}}^S,\emph{\textbf{x}}^C}-2\
\hspace*{-0.1cm}b\
\hspace*{-0.1cm}\widehat{\sigma}_{\emph{\textbf{x}}^S,\emph{\textbf{y}}} -2\
\hspace*{-0.1cm}c\
\hspace*{-0.1cm}\widehat{\sigma}_{\emph{\textbf{x}}^C,\emph{\textbf{y}}}.
\end{equation} 

\noindent If the global minimum  of the function $g$ is so that
 $(b^*,c^*)^t \notin \Gamma_G$, then the local minimum
of $g$ over $\Gamma_G$ is unique, and it is located on the boundary of
$\Gamma_G$. The boundary of $\Gamma_G$, denoted by $\mathrm{fr}(\Gamma_G)$, verifies that
\begin{equation}\label{frontera-curvas}
\mathrm{fr}(\Gamma_G)=L_1 \cup L_2 \cup L_3\ ,
\end{equation}
\noindent where $L_i$, $i=1,2,3$ are the following sets:

\begin{itemize}
\item $L_1=\left\{(0,c)\,|\,0\leq c \leq r_0=\min_{i=1,\ldots,n} \frac{\mathrm{spr}
\emph{\textbf{y}}_i}{|\mathrm{mid}
\emph{\textbf{x}}_i|}\right\}$.

\item $L_3=\left\{(b,0)\,|\, 0\leq b \leq s_0=\min_{i=1,\ldots,n} \frac{\mathrm{spr}
\emph{\textbf{y}}_i}{\mathrm{spr}\emph{\textbf{x}}_i} \right\}$.

\item $L_2=\left\{(b,\min_{k=1 \ldots n} \{-u_k b + v_k\}) \,|\,0\leq b \leq s_0 \right\}$, with \\
$$\displaystyle u_k=\frac{\mathrm{spr} \emph{\textbf{x}}_k}{|\mathrm{mid} \emph{\textbf{x}}_k|} \ \mbox{ and } \
v_k=\frac{\mathrm{spr} \emph{\textbf{y}}_k}{|\mathrm{mid} \emph{\textbf{x}}_k|} \textrm{ for all } k=1,\ldots,n.$$
\end{itemize}

The set $L_2$ is composed on several straight segments from some of the
straight lines $\{l_k: c= - u_k b + v_k\}_{k=1}^n$. If $|\mathrm{mid} \emph{\textbf{x}}_k| =0$
for any $k\in \{1,\ldots, n\}$, then the corresponding straight line is
$b=\mathrm{spr} \emph{\textbf{y}}_k/\mathrm{spr} \emph{\textbf{x}}_k$ for $\mathrm{spr} \emph{\textbf{x}}_k\neq 0$. Thus,
it is a vertical line, which could take part in   
 $L_2$ only if $\mathrm{spr} \emph{\textbf{y}}_k/\mathrm{spr} \emph{\textbf{x}}_k=s_0$. Moreover, if
$\mathrm{spr} \emph{\textbf{x}}_k = 0$ too, then the sample interval $\emph{\textbf{x}}_k$ is reduced to the
real value $\emph{\textbf{x}}_k=0$, so it does not take part in the construction of $\Gamma_G$. In Figure 1 the feasible set
and its boundary in a practical example are illustrated graphically. The sample
data corresponds to a real-life example (see Section~\ref{subsection-example}). 

\begin{figure}
\begin{center}
\includegraphics[height=7cm,width=7cm]{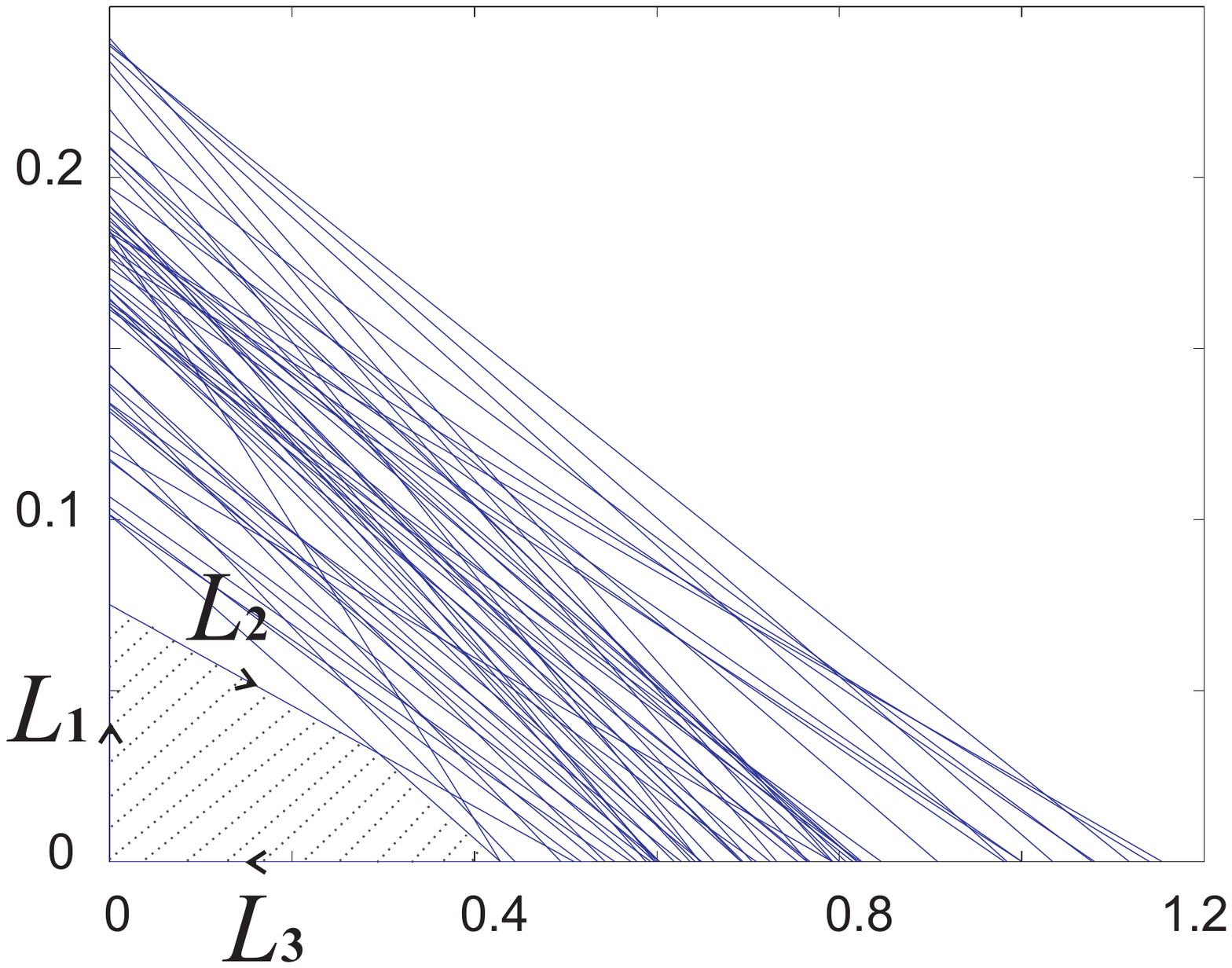}
\vspace*{-0.5cm}
\caption{\small $\Gamma_G$ for the sample data in Table   
\ref{table-hospital}}\label{fig_feasible_hospital}
\end{center}
\end{figure}

In order to find the exact solution of 
$\min_{(b,c)\in \Gamma_G} g(b,c)$
the global minimum of  $g$  should be computed and, if needed, the
local minimum over $L_i, i=1,2,3$.

The asymptotic time complexity of the 
algorithm is $0(nt)$, where $t$ is the number of lines in $\{l^k\}_{k=1}^{n}$ taking part in 
$\mathrm{fr}(\Gamma_G)$. The straight lines in $\{l_k: c= - u_k b + v_k: k\neq
(v),(h)\}_{k=1}^n$ such that $-u_k b_{(v,h)} + v_k > c_{(v,h)}$ do not take part on
the construction of $\mathrm{fr}(\Gamma_G)$. Thus, they can be ignored from Step 5 to the end of the
algorithm.  However, for practical
examples
with moderate sample sizes $n$, this reduction will result in a negligible improvement on the
computational efficiency of the algorithm.

 \vspace*{0.5cm}

\section{The multiple basic linear regression
model}\label{subsection-basic-multiple}

 Let \emph{\textbf{y}} be a response random interval and let 
$\emph{\textbf{x}}_1,\emph{\textbf{x}}_2,\ldots,\emph{\textbf{x}}_k$ be $k$ explanatory random
intervals.
The multiple basic linear regression model (MBLRM) extending (\ref{model-basic})  is
formalized as:
\begin{equation}\label{model-basic-multiple}
\emph{\textbf{y}}=x^t b+\boldsymbol{\varepsilon}
\end{equation}
\noindent being
$x=(\emph{\textbf{x}}_1,\emph{\textbf{x}}_2,\ldots,\emph{\textbf{x}}_k)^t$,
$b=(b_1,b_2,\ldots,b_k)^t \in \mathbb{R}^k$ and $\boldsymbol{\varepsilon}$ an random
interval-valued error such that $E(\boldsymbol{\varepsilon}|x)=\Delta \in \mathcal{K
_{\mathrm{c}}}(\mathbb{R})$.
The associated regression function
is
$E(\emph{\textbf{y}}|\emph{\textbf{x}}_1=x_1,\ldots,\emph{\textbf{x}}
_k=x_k)=x^tb+\Delta$.  Thereafter, the second-order moments of the random intervals involved
in the linear model (\ref{model-basic-multiple}) are assumed to be finite, and
the variances strictly positive. If the  mids and spreads of the explanatory intervals are not
degenerated, then (\ref{model-basic-multiple}) is unique. The following separate models  are transferred:
\begin{equation}\label{separate-basic-multiple}
\left.\begin{array}{c}
\mathrm{mid}\emph{\textbf{y}}=\mathrm{mid}(x^t)\,b+\mathrm{mid}\,\boldsymbol{\varepsilon}, \textrm{ and } \\
\mathrm{spr}\emph{\textbf{y}}=\mathrm{spr}(x^t)\,|b|+\mathrm{spr}\,\boldsymbol{\varepsilon}.  \end{array}\right.
\end{equation}
\noindent The {\it mid} variables relates through a
standard (real-valued) multiple linear model, but this is not the case for the  spreads, due to the non-negative restrictions.

Let $\{\left(\emph{\textbf{y}}_j,\emph{\textbf{x}}_{i,j}\right) : i=1,\ldots,k, j
=1,\ldots,n\}$ be a simple random sample of size $n$ obtained from $\emph{\textbf{y}}$ and $x=(\emph{\textbf{x}}_1,\ldots,
\emph{\textbf{x}}_k)$ verifying (\ref{model-basic-multiple}). Then, 
\begin{equation}\label{sample-equation-basic-multiple}
y=X b+\varepsilon,
\end{equation}
where $y=(\emph{\textbf{y}}_1,\ldots,\emph{\textbf{y}}_n)^t$, $X$ is the $(n
\times k)$-interval-valued matrix such that $X_{j,i}= \emph{\textbf{x}}_{i,j}$,
and $\varepsilon=(\boldsymbol{\varepsilon}_1,\ldots,\boldsymbol{\varepsilon}_n)^t$ fulfils
$E[\varepsilon|X]= 1^n \Delta$, $1^n$ denoting the ones' vector in
$\mathbb{R}^n$. The LS estimation  consists in finding
$\widehat{b}$ and $\widehat{\Delta}$ minimizing the  objective function:
\begin{equation}\label{min1-basic-multiple}
\min_{d\in\mathbb{R}^k, C \in \mathcal{K _{\mathrm{c}}}(\mathbb{R})}
d_{\theta}^2(y, X d +1^n C)\ ,
\end{equation}
constrained to the existence of the residuals
$\widehat{\varepsilon}=y-_{H}X\widehat{b}$. If $\widehat{b}\in \Gamma=\{a \in \mathbb{R}^k: y-_{H}X
a\in \mathcal{K _{\mathrm{c}}}(\mathbb{R})^k\}$, then the optimum value  over $C$  is attained at
\begin{equation}\label{delta-basic-multiple}
\widehat{\Delta}=\overline{\emph{\textbf{y}}} -_H \overline{x^t}\widehat{b}.
\end{equation}

Extending directly the estimation method in \cite{GRetal:07} would lead to a computationally
unfeasible combinatorial problem. For that, a non-optimal stepwise algorithm  has been  proposed.
However, that may be offset by estimating separately the absolute
value of $\widehat{b}$ and its sign. Note that $\widehat{b}=|\widehat{b}| \circ
sign(\widehat{b})$, and from (\ref{separate-basic-multiple}),
$sign(\widehat{b})$ is only determined by the sign of the relationship between the mid-points. Then,
$sign(\widehat{b})_i=sign(\widehat{Cov}(\textrm{mid}\emph{\textbf{y}},
\textrm{mid}x_i))$ and $|\widehat{b}|$ can be obtained as the solution of
\begin{equation}\label{min2-basic-multiple}
\min_{a\in\Gamma, a\ge 0}d_{\theta}^2(y, X( a \circ sign(\widehat{b})) +1^n
\widehat{\Delta}) .
\end{equation}

\noindent The feasible set $\Gamma'= \Gamma \cap (\mathbb{R}^k)^{+}$ in
(\ref{min2-basic-multiple}) can be expressed as
\begin{equation}\label{gammaprima}
\Gamma'= \{d\in (\mathbb{R}^k)^{+} : (\textrm{spr} X) d \le \textrm{spr} y\}.
\end{equation}

\noindent  A more operative expression for the objective
function in (\ref{min2-basic-multiple}) is:
\begin{eqnarray}\label{min3}
\nonumber \hspace*{-0.3cm} d_{\theta}^2(y, X (a \circ
sign(\widehat{b}))\hspace*{-0.1cm} + \hspace*{-0.1cm}1^n
\widehat{\Delta})&\hspace*{-0.3cm}=\hspace*{-0.3cm}&(v_m-F_m'\,a)^t(v_m-F_m'\, a)\\
\hspace*{-0.3cm} & \hspace*{-0.3cm} \hspace*{-0.3cm}& \hspace*{-0.5cm}+\theta
(v_s-F_s' \,a)^t (v_s-F_s'\, a),
\end{eqnarray}
\noindent where $v_m=\mathrm{mid}y-\overline{\mathrm{mid}\emph{\textbf{y}}}1^n\in \mathbb{R}^n$,
$v_s=\mathrm{spr}y-\overline{\mathrm{spr}\emph{\textbf{y}}}1^n\in \mathbb{R}^n$,
$F_m=(\mathrm{mid}X-1^n(\overline{\mathrm{mid}x^t}))$,
$F_s=\mathrm{spr}X-1^n(\overline{\mathrm{spr}x^t})$, $F_m'=F_m \textrm{diag}(sign(\widehat{b})_1,\ldots, 
sign(\widehat{b})_k)\in\mathbb{R}^{n\times k}$, 
$F_s'=F_s \textrm{diag}(sign(\widehat{b})_1,\ldots, 
sign(\widehat{b})_k)\in \mathbb{R}^{n\times k}$ and $a\in\mathbb{R}^k$. Since the optimization problem consists in
minimizing a quadratic expression with inequality linear constraints, Karush-Kuhn-Tucker (KKT) conditions
guarantee the
existence of solution and it can be found by using a standard software.

\section{The multiple flexible linear regression
model}\label{subsection-MG-multiple}

From (\ref{model-Mg}), a
multiple flexible linear regression model (MFLRM) can be defined as:
\begin{equation}\label{model-MG-multiple}
\emph{\textbf{y}}=\mathrm{mid}\,x^t\,[1\pm 0]\,b^1+\mathrm{spr}\,x^t\,[0\pm
1]\,b^2+\mathrm{mid}\,x^t\,[0\pm 1]\,b^3+\mathrm{spr}\,x^t\,[1\pm
0]\,b^4+\boldsymbol\varepsilon \ , 
\end{equation}
where $b^1,b^2,b^3,b^4\in \mathbb{R}^k$ and $E(\boldsymbol\varepsilon|x^t)=\Delta \in
\mathcal{K _{\mathrm{c}}}(\mathbb{R})$. Equivalently (\ref{model-MG-multiple}) can be written  as:
\begin{equation}\label{MFLRM-intervals}
\emph{\textbf{y}}=x^M\,b^1+x^S\,b^2+x^C\,b^3+x^R\,b^4+\boldsymbol{\varepsilon}  ,
\end{equation}
\noindent or, in matrix notation, as:
\begin{equation}\label{MFLRM-matrix}
\emph{\textbf{y}}=X^{Bl}\,B+\boldsymbol{\varepsilon} \ ,
\end{equation}
where $X^{Bl}=(x^M|x^S|x^C|x^R)\in\mathcal{K_{\mathrm{c}}}(\mathbb{R})^{1\times
4k}$ and $B=((b^1)^t|(b^2)^t|(b^3)^t|(b^4)^t)^t\in \mathbb{R}^{4k \times 1}$.
The values $b^2$ and $b^3$ can be assumed to be non-negative without loss of
generality since $x^S=-x^S$ and $x^C=-x^C$.

The separate linear relationships for the {\it mid} and {\it spr}
components of the intervals transferred from (\ref{model-MG-multiple}) are 
\begin{equation}
\mathrm{mid}\,\emph{\textbf{y}}=\mathrm{mid}\,(x^t)\,b^1+\mathrm{spr}\,(x^t)\,
b^4+\mathrm{mid}\,\boldsymbol{\varepsilon}\ \mathrm{, and}
\end{equation}
\begin{equation}
\mathrm{spr}\,\emph{\textbf{y}}=\mathrm{spr}\,(x^t)\,b^2+|\mathrm{mid}\,(x^t)|\,
b^3+\mathrm{spr}\,\boldsymbol{\varepsilon}.
\end{equation}

Let $\{\left(\emph{\textbf{y}}_j,\emph{\textbf{x}}_{i,j}\right) :
i=1,\ldots,k, j =1,\ldots,n\}$ be a simple random
sample obtained from the random intervals $(\emph{\textbf{y}},
\emph{\textbf{x}}_{1},\ldots,\emph{\textbf{x}}_{k})$ verifying 
(\ref{model-MG-multiple}). Then,
$$y=X^M\,b^1+X^S\,b^2+X^C\,b^3+X^R\,b^4+\varepsilon\ ,$$
\noindent where $y=(\emph{\textbf{y}}_1,\ldots,\emph{\textbf{y}}_n)$, 
$\varepsilon=(\boldsymbol\varepsilon_{1},\ldots,\boldsymbol\varepsilon_{n})$ such that
$E(\varepsilon
|x)=1^n \Delta$, $(X^M)_{i,j}=\mathrm{mid}\emph{\textbf{x}}_i [1\pm 0]$ and
$X^S,X^C$ and $X^R$ are analogously defined. It can be equivalently expressed in a
matrix form as 
\begin{equation}\label{MFLRM-matricial-sample}
y=X^{ebl}B+\varepsilon \ ,
\end{equation}
\noindent where $X^{ebl}=(X^M|X^S|X^C|X^R)\in \mathcal{K
_{\mathrm{c}}}(\mathbb{R})^{n\times 4k}$ and $B$ as in (\ref{MFLRM-matrix}). 

The LS estimation  searches for $\widehat{B}$ and
$\widehat{\Delta}$ minimizing 
$\displaystyle d_{\theta}^2(y, X^{ebl}A+1^n C)$ for $A\in
\mathbb{R}^{4k\times 1}$ and $C\in \mathcal{K}_c(\mathbb{R})$. The constraints to assure the
existence of the residuals are: 
\begin{equation}\label{conditions}
\mathrm{spr}\,X\,\widehat{b^2}+|\mathrm{mid}\,X|\,\widehat{b^3}\leq
\mathrm{spr}\,y .
\end{equation}

The estimation of $B$ and $\Delta$ can be solved separately. If $\widehat{B}$ verifies
(\ref{conditions}), then the minimum value of $d_{\theta}^2(y,
X^{ebl}\widehat{B}+1^n C)$ over $C\in \mathcal{K}_c(\mathbb{R})$ is attained at
$\widehat{\Delta}= \overline{\emph{\textbf{y}}}-_H
\overline{X^{Bl}}\widehat{B}$.  The objective function can then be written as $$d_{\theta}^2(y,
X^{ebl}A+1^n \widehat{\Delta})=
(v_m-F_mA_m)^t(v_m-F_mA_m)+\theta\,(v_s-F_sA_s)^t(v_s-F_sA_s)\ ,$$

\noindent where $v_m, v_s \in \mathbb{R}^n, F_m,\, F_s\in\mathbb{R}^{n\times 2k}$ are defined as in (\ref{min3}), $A_m=((a^2)^t|(a^3)^t)^t\in \mathbb{R}^{2k}$ are the coefficients affecting the mids and $A_s=((a^1)^t|(a^4)^t)\in\mathbb{R}^{2k}$ the coefficients affecting the spreads, with $a^l\in \mathbb{R}^k$, $l=1,\ldots,4$. 

Therefore, the computation of the LS estimator $\widehat{B}$ for the regression parameter $B$ in (\ref{MFLRM-matrix}) is solved through the constrained optimization problem by KKT conditions:  
\begin{equation}\label{minim-MFLRM}
\min_{A_s\in\mathbb{R}^{2k}, A_m\in \Gamma_2}(v_m-F_m A_m)^t(v_m-F_m
A_m)+\theta\,(v_s-F_s A_s)^t(v_s-F_s A_s)\ , 
\end{equation}
\noindent with 
\begin{equation}\label{Gamma2}
\Gamma_2=\{(a^2,a^3)\in [0,\infty)^k\times[0,\infty)^k :
\mathrm{spr}\,X\,a^2+|\mathrm{mid}\,X|\,a^3\leq \mathrm{spr}\,y\}.
\end{equation}

Note that the extension of the linear regression model M developed in
\cite{Blancoetal:11} to the multiple case is directly achieved from
(\ref{model-MG-multiple}), taking $b^3=(0,\ldots,0)$ and $b^4=(0,\ldots,0)$.

\section{Empirical results}\label{section-empirical}

\subsection{Application to a real-life example}\label{subsection-example}

A real-life example concerning the relationship between the daily
fluctuations of the systolic and diastolic blood
pressures and the pulse rate over a sample of patients in the
Hospital {\it Valle del Nal{\'o}n}, in Spain, is considered (previously explored in \cite{Blancoetal:11, Giletal:02, GRetal:07}). The metric $d_{1/3}$ is employed, and the optimization
algorithm {\it quadprog} to solve the estimation of the multiple models (\ref{model-basic-multiple}) and (\ref{model-MG-multiple}) is
run.
 
Let $\emph{\textbf{y}}$, $\emph{\textbf{x}}_1$ and  $\emph{\textbf{x}}_2$  the {\it fluctuation of
the diastolic blood
pressure of a patient over a day}, the {\it fluctuation of the systolic
blood pressure over the same day}, and the {\it pulse range
variation over the same day}, respectively. Data in Table \ref{table-hospital} correspond to a
sample data of $59$ patients from
$(\emph{\textbf{y}},\emph{\textbf{x}}_1,\emph{\textbf{x}}_2)$.

From the sample data provided in Table
\ref{table-hospital}, the estimated model $\mathrm{M}_G$ for $\emph{\textbf{y}}$
and  $\emph{\textbf{x}}_1$ is 

\begin{equation}\label{example-MG}
\widehat{\emph{\textbf{y}}} = 0.5383\emph{\textbf{x}}_1^M +
0.2641\emph{\textbf{x}}_1^S- 0.4412\emph{\textbf{x}}_1^R + [4.249, 35.254].
\end{equation}

\noindent The value of determination coefficient $R^2$ (defined as the proportion of
explained variability) associated  with this estimated
model is $0.6857$.
 
 The estimated model
(\ref{model-basic-multiple}) for  $\emph{\textbf{y}}$ and $(\emph{\textbf{x}}_1,\emph{\textbf{x}}_2)$ from the data set in Table \ref{table-hospital} has
the expression:  

\begin{equation}\label{example-MBLRM}
\widehat{\emph{\textbf{y}}} = 0.4094\emph{\textbf{x}}_1 +
0.0463\emph{\textbf{x}}_2 + [10.3630, 29.5168].
\end{equation}

\noindent The value of the determination coefficient is in this case
$R^2=0.4221$.

 The linear relationship between ${\emph{\textbf{y}}}$ and
$(\emph{\textbf{x}}_1,\emph{\textbf{x}}_2)$ can be also estimated more naturally  by means of
the MFLRM. The estimation of the model
(\ref{model-MG-multiple}) leads to the expression: 
\begin{eqnarray}\label{example-MFLRM}
\nonumber\widehat{\emph{\textbf{y}}} &\hspace*{-0.3cm}
=\hspace*{-0.3cm}&0.5435\,\emph{\textbf{x}}_1^M+ 0.0190\,\emph{\textbf{x}}_2^M+
0.2588\,\emph{\textbf{x}}_1^S+
0.1685\,\emph{\textbf{x}}_2^S-0.4446\,\emph{\textbf{x}}_1^R\\
 & & +0.1113\,\emph{
\textbf{x}}_2^R + [3.2032, 27.8373]\ , 
\end{eqnarray}

\noindent with $R^2=0.7922$.
\medskip

The highest value of $R^2$ is
achieved for (\ref{example-MFLRM}), which agrees with the fact that MFLRM is the most
flexible regression among the linear models that have been developed. The difference in
the $R^2$ between this multiple model and the simple one in (\ref{example-MG})
is due to the inclusion of the pulse rate variable $\emph{\textbf{x}}_2$ in
the prediction of $\emph{\textbf{y}}$. However, this
difference is not large, which indicates that the pulse rate has low explanatory power. The smallest
value of $R^2$ corresponds to
(\ref{example-MBLRM}). It indicates that the multiple basic model is too restrictive to relate
these physical magnitudes.

\subsection{Simulation results}\label{subsection-simulations}

The empirical performance of the regression estimates for each linear model is investigated
 by means of some simulations. Three independent random intervals $\emph{\textbf{x}}_1,
\emph{\textbf{x}}_2, \emph{\textbf{x}}_3$ and an interval error $\boldsymbol\varepsilon$
will be considered.  Let $\mathrm{mid}\,\emph{\textbf{x}}_1\sim
\mathcal{N}(1,2)$, $\mathrm{spr}\,\emph{\textbf{x}}_1\sim \mathcal{U}(0,10)$,
$\mathrm{mid}\,\emph{\textbf{x}}_2\sim \mathcal{N}(2,1)$,
$\mathrm{spr}\,\emph{\textbf{x}}_2\sim \mathcal{X}_4^2$,
$\mathrm{mid}\,\emph{\textbf{x}}_{\,3}\sim \mathcal{N}(1,3)$,
$\mathrm{spr}\,\emph{\textbf{x}}_{ \,3}\sim \mathcal{U}(0,5)$,
$\mathrm{mid}\,\boldsymbol\varepsilon\sim \mathcal{N}(0,1)$ and
$\mathrm{spr}\,\boldsymbol\varepsilon\sim \mathcal{X}_1^2$. Different linear expressions
with the investigated structures will be considered.

\begin{itemize}
\item Model $M_1$: According to the multiple basic linear model
presented in (\ref{model-basic-multiple}), $\emph{\textbf{y}}$ it is defined by
the expression:
\begin{equation}\label{sim1}
\emph{\textbf{y}}=2\emph{\textbf{x}}_1-5\emph{\textbf{x}}_2-\emph{\textbf{x}}_{3}+\varepsilon.
\end{equation}

\item Model $M_2$: A simple linear relationship in terms of the model
$\mathrm{M}_G$ in (\ref{model-Mg}) is defined by considering only
$\emph{\textbf{x}}_1$ as independent interval for modelling $\emph{\textbf{y}}$
through the expression:
\begin{equation}\label{sim2}
\emph{\textbf{y}}=-2\emph{\textbf{x}}_1^M+2\emph{\textbf{x}}_1^S+\emph{\textbf{x}}_1^C+0.5\emph{
\textbf{x}}_1^R+
\varepsilon.
\end{equation}

\item Model $M_3$: A multiple flexible linear regression model following
(\ref{model-MG-multiple}) is defined as:
\begin{eqnarray}\label{sim3}
\nonumber
\emph{\textbf{y}}&=&-2\emph{\textbf{x}}_1^M+5\emph{\textbf{x}}_2^M-\emph{\textbf
{x}}_3^M+2\emph{\textbf{x}}_1^S+2\emph{\textbf{x}}_2^S+\emph{\textbf{x}}
_3^S+\emph{\textbf{x}}_1^C+\emph{\textbf{x}}_2^C+3\emph{\textbf{x}}_3^C\\
& &
+0.5\emph{\textbf{x}}_1^R+\emph{\textbf{x}}_2^R-3\emph{\textbf{x}}
_3^R+\varepsilon.
\end{eqnarray}

\end{itemize}

From each linear model $l=10,000$ random samples has been generated for different sample
sizes $n$. The estimates of the regression parameters have been computed for each iteration.
Table \ref{table-sim}
shows the estimated mean value and MSE of the LS
estimators (denoted globally by $\widehat{\nu}$) computed from the $l$
iterations. The mean values of
the estimates are always closer to the corresponding regression parameters as the
sample size $n$ increases, which empirically shows the asymptotic unbiasedness
of the estimators. Moreover, the values for the estimated MSE tend to zero as $n$
increases.

 In Figure 2 the box-plots of the $l$ estimates of the model $M_1$
are presented for $n=30$ (left-side plot)  and  $n=100$
(right-side plot) sample observations. All the cases the boxes reduce their
width around the true value of the corresponding parameter on the population
linear model as the sample size $n$ increases, which illustrates the consistency. Analogous
conclusions are obtained
for the models $M_2$ and $M_3$ in Figures 3 and 4,
respectively.

\begin{figure}[!h]
\begin{center}
\begin{minipage}{7cm}
\hspace*{-2cm}
\includegraphics[height=12cm,width=10cm]{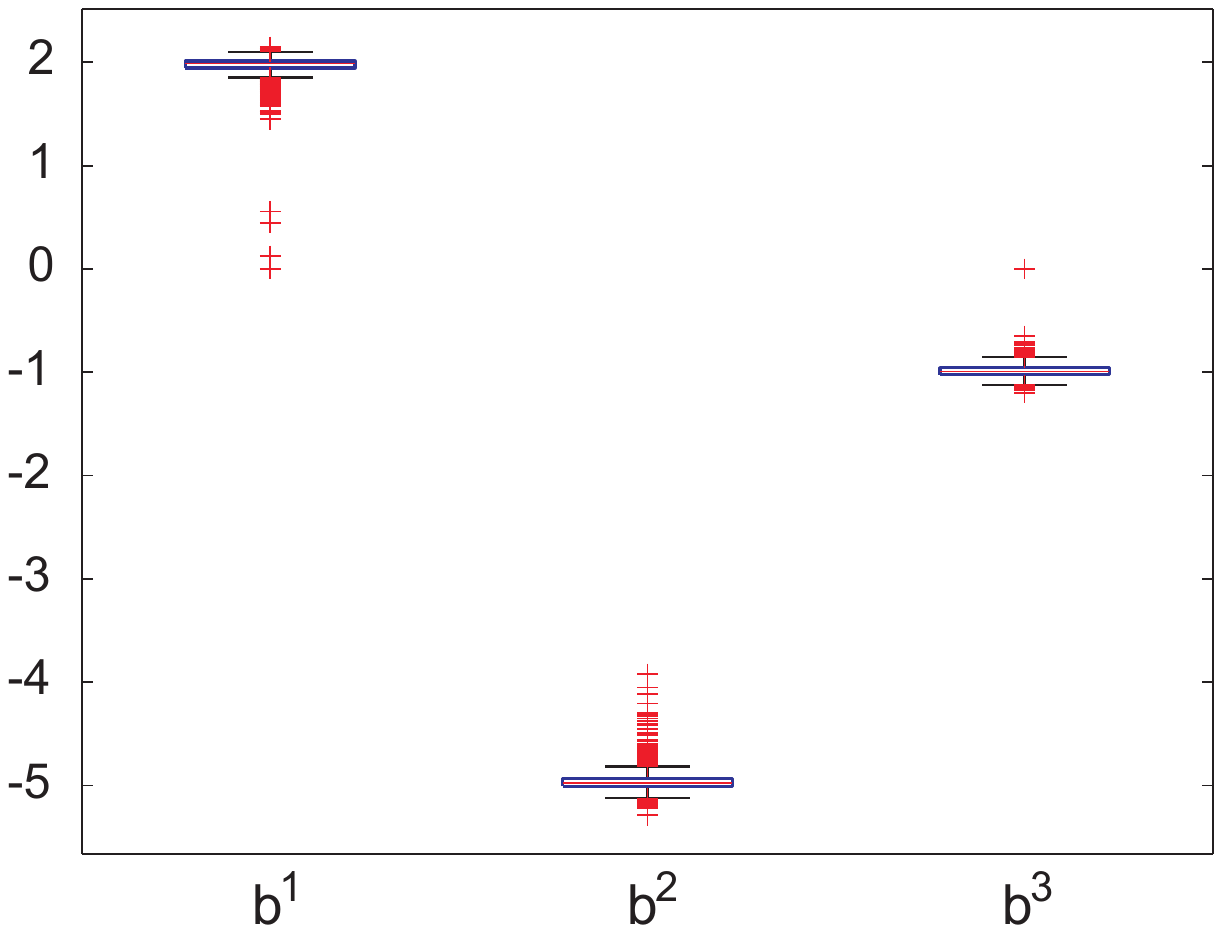}
\end{minipage}
\begin{minipage}{6cm}
\hspace*{-2cm}
\includegraphics[height=12cm,width=10cm]{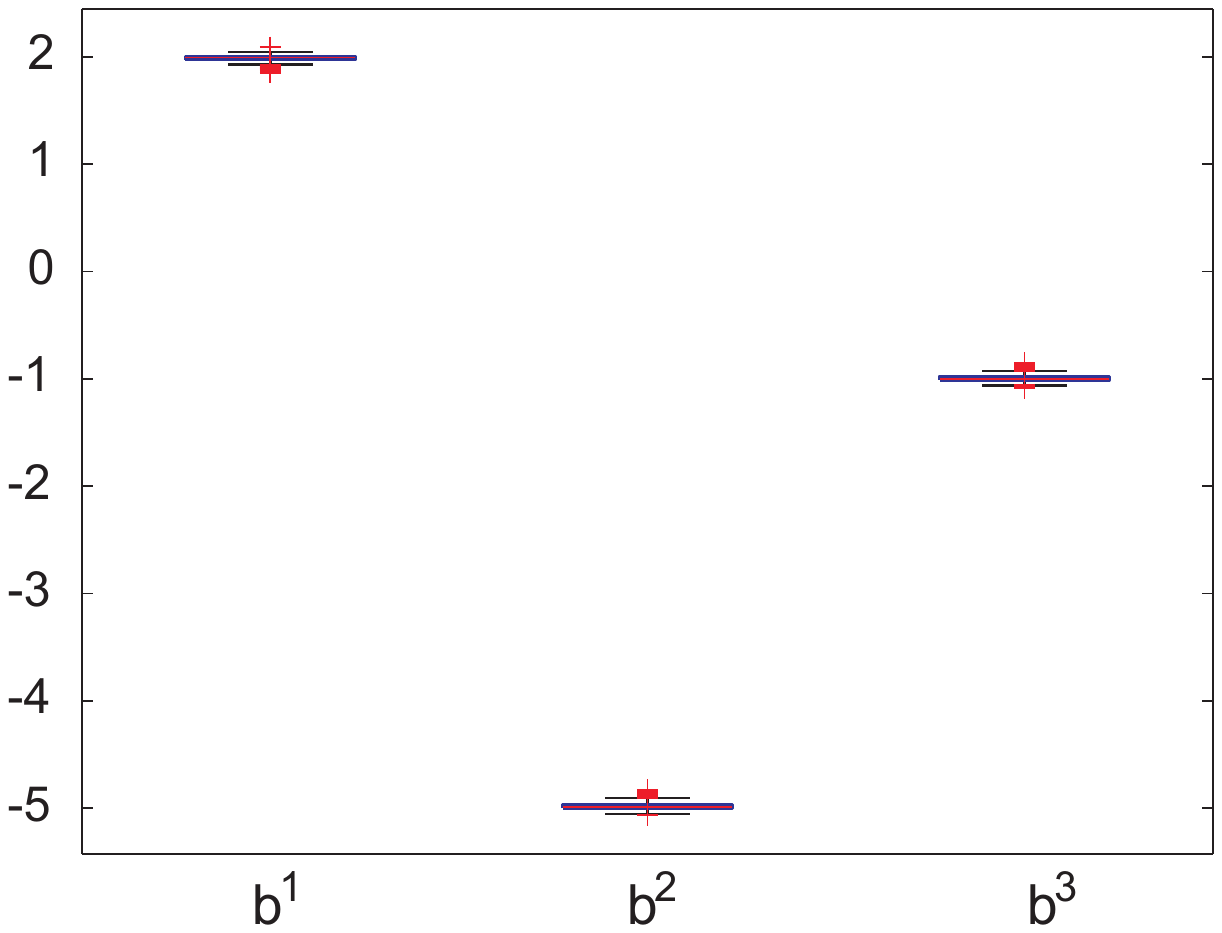}
\end{minipage}
\vspace*{-4cm}
\caption{\small Box plot of the LS estimators for model $M_1$, $n$=30 (left);
$n$=100 (right)}\label{figures1}
\end{center}
\end{figure}

\begin{figure}[!h]
\begin{center}
\vspace*{-4cm}
\begin{minipage}{7cm}
\hspace*{-2cm}
\includegraphics[height=12cm,width=10cm]{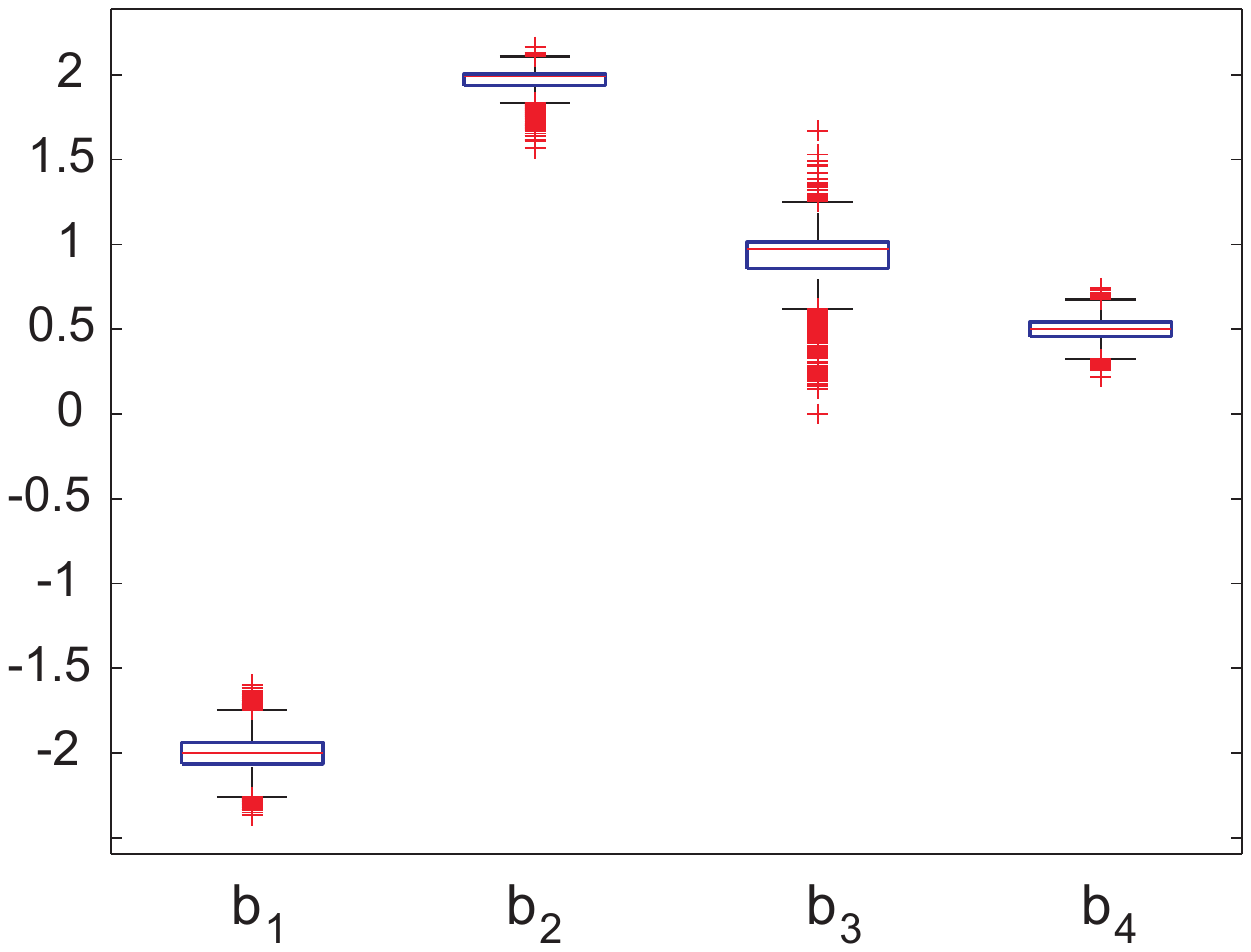}
\end{minipage}
\begin{minipage}{6cm}
\hspace*{-2cm}
\includegraphics[height=12cm,width=10cm]{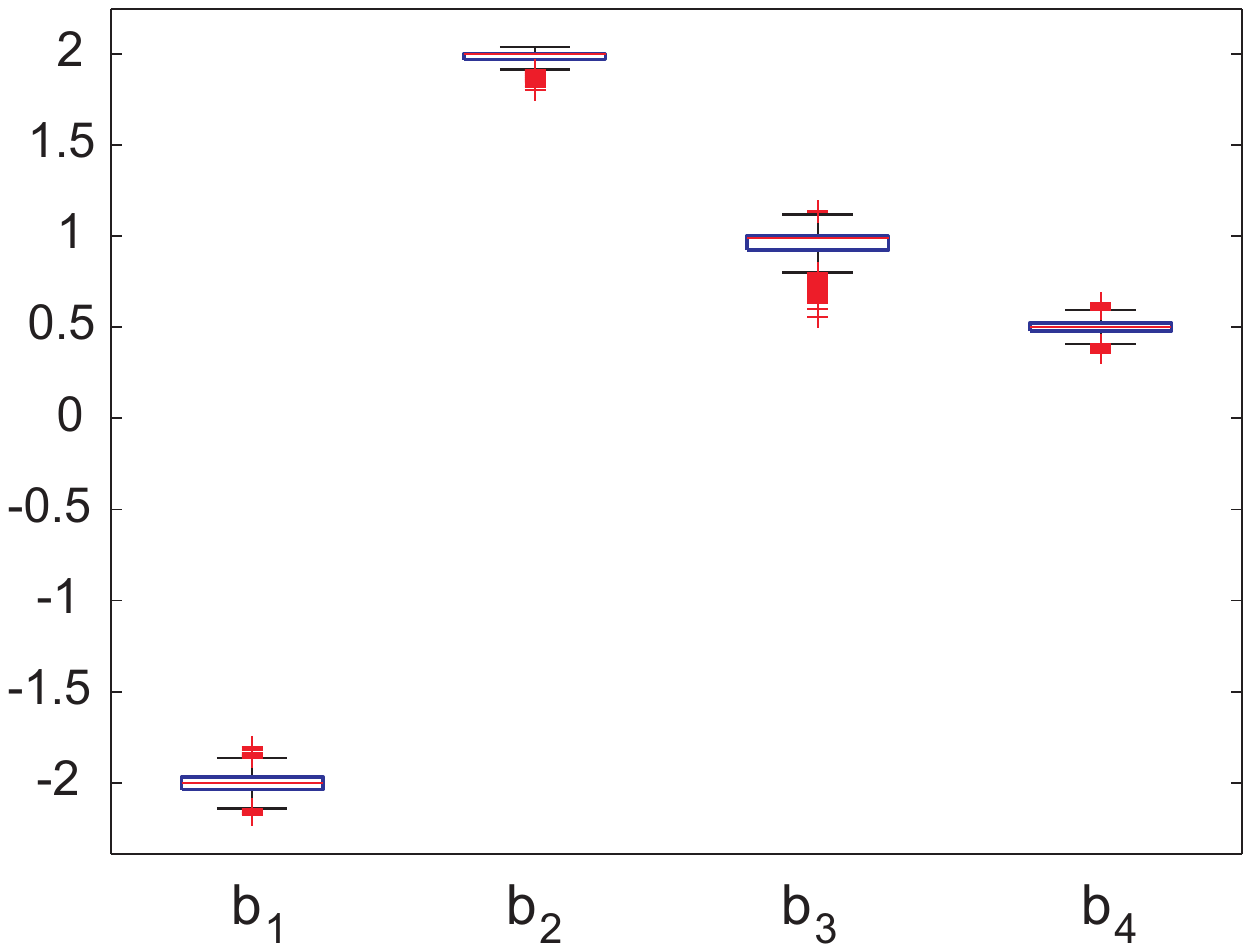}
\end{minipage}
\vspace*{-4cm}
\caption{\small Box plot of the LS estimators for model $M_2$, $n$=30 (left);
$n$=100 (right)}\label{figures2}
\end{center}
\end{figure}
\begin{figure}[!h]
\begin{center}
\vspace*{-4cm}
\begin{minipage}{7cm}
\hspace*{-2cm}
\includegraphics[height=12cm,width=10cm]{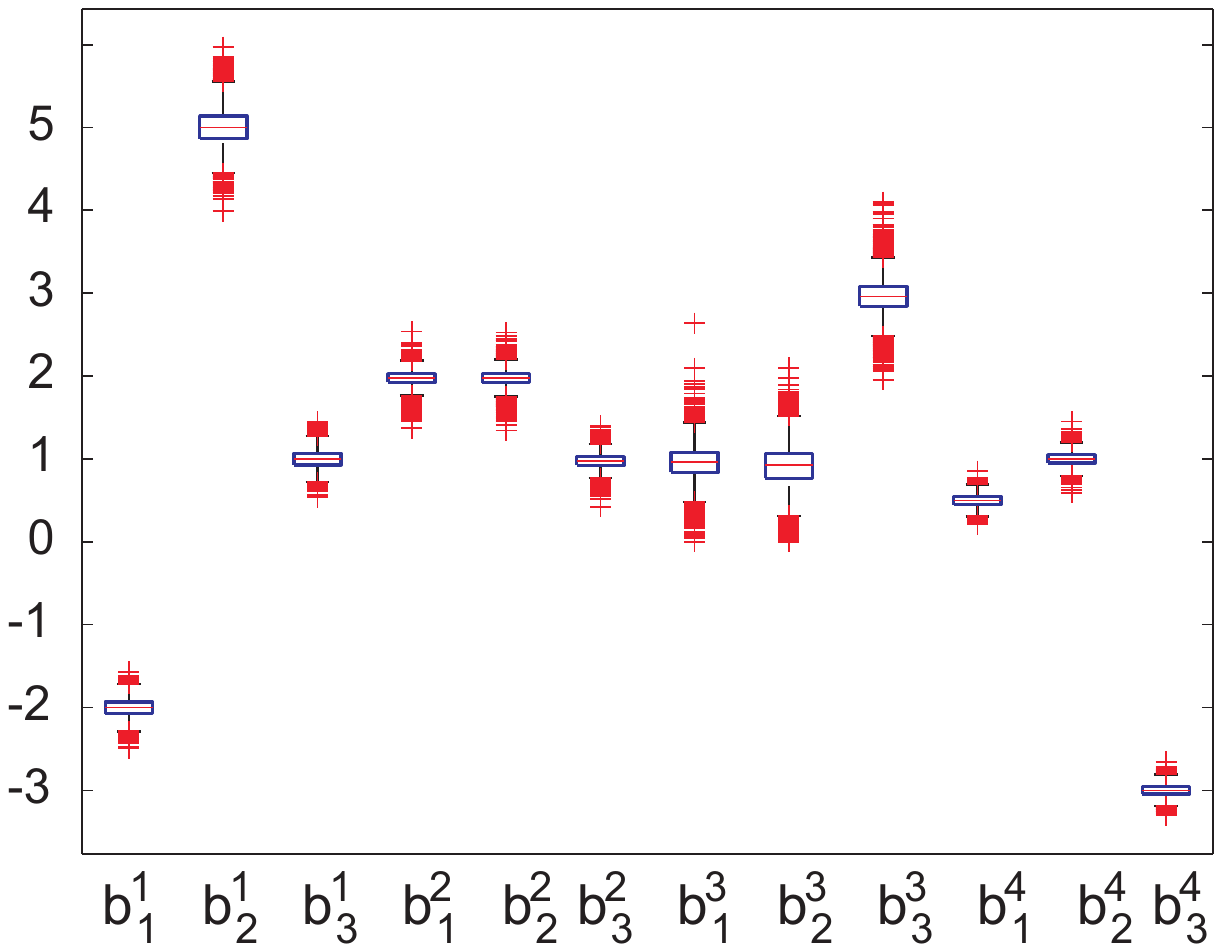}
\end{minipage}
\begin{minipage}{6cm}
\hspace*{-2cm}
\includegraphics[height=12cm,width=10cm]{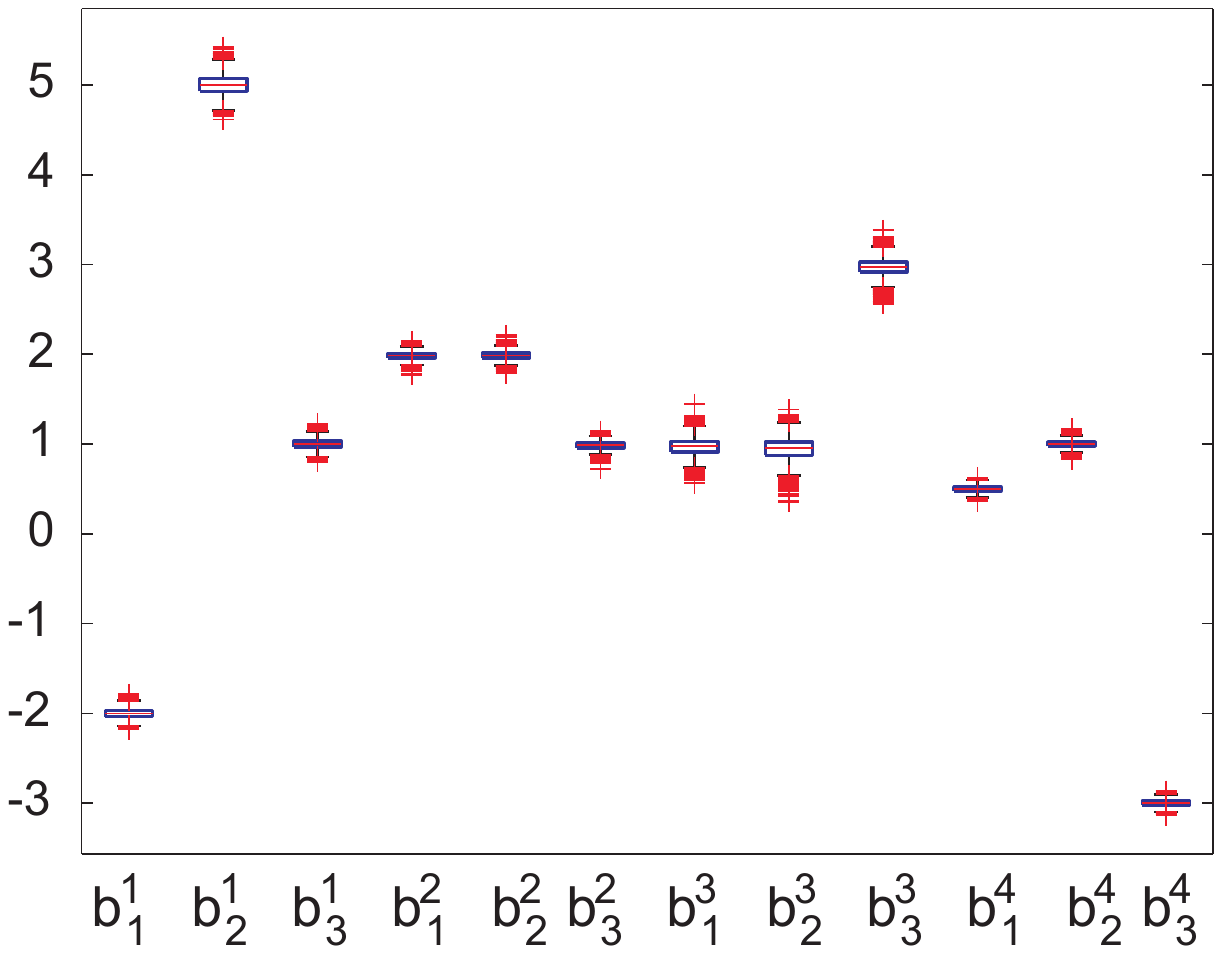}
\end{minipage}
\vspace*{-4cm}
\caption{\small Box plot of the LS estimators for model $M_3$, $n$=30 (left); 
$n$=100 (right)}\label{figures3}
\end{center}
\vspace*{-0.5cm}
\end{figure}

\section{Conclusions}\label{section-conclusions}

Previous linear regression models for interval data based on
set arithmetic have been extended.  
In all cases the search of the LS estimators involves minimization problems with constraints. 
The constraints are necessary to assure the existence of the residuals and thus, the coherency 
of the estimated model with the population one. 

A very flexible simple model based on the canonical decomposition and allowing
for cross-relationships between mid-points and spreads has been introduced. 
An algorithm to find the exact LS-estimates has been developed. This model has been extended
to the multiple case. The LS exact algorithm strongly relies on the geometry of the feasible set 
and it cannot be generalized in a simple way. However, the LS estimates can be found by applying 
the KKT conditions. The extension of the basic simple model in \cite{GRetal:07}, which is not based on the canonical 
decomposition, requires a different approach, but the solutions can also be found by applying 
the KKT conditions.

The empirical validity of the estimation process for all the models has been
shown by means of  simulations. However, further
theoretical studies of the main properties of the regression estimators, as the
bias, the consistency or the asymptotic distributions should be pursued.

\section*{Acknowledgements}
 
The research in this paper has
been partially supported by the Spanish Ministry of Science and
Innovation Grant MTM2009-09440-C02-01. It has also
benefited from short-term scientific missions associated with the COST
Action IC0702. 








\newpage
\setcounter{table}{0}
\begin{table}[!h]
\caption{\small Sample data for the daily blood pressures and the pulse rate ranges of
$59$ patients}\label{table-hospital}
\vspace{0.5cm}
\begin{center}
{\small \begin{tabular}{ c  c  c | c  c  c | c  c  c }
\hline
$\emph{\textbf{y}}$ & $\emph{\textbf{x}}_1$ & $\emph{\textbf{x}}_2$ &
$\emph{\textbf{y}}$ &  $\emph{\textbf{x}}_1$  & $\emph{\textbf{x}}_2$  &
$\emph{\textbf{y}}$  & $\emph{\textbf{x}}_1$ & $\emph{\textbf{x}}_2$\\
\hline
63-102 &118-173 &58-90  &47-93 &119-212 &52-78 & 71-118 &104-161 &47-68\\
 73-105&122-178 &55-84 & 58-113 &131-186 &32-114 & 74-125 &127-189 &61-101 \\
62-118 &105-157 &61-110 &52-112 &113-213 &65-92 &59-94  &120-179 &62-89\\
69-133 &141-205 &38-66 &48-116 &101-194 &63-119 &53-109 &99-169 &48-73 \\
60-119 &109-174 &51-95 &60-98 &126-191 &59-98 & 76-125&128-210 &49-78 \\
55-121 &99-201 &59-87 &47-104 &94-145 &43-67 &37-94 &88-221 &49-82 \\
88-130 &148-201 &55-102 & 55-85 &113-183 &48-77&52-96  &111-192 &64-107 \\
56-121 &94-176 &56-133 & 74-133&116-201 &54-84 & 50-94&102-156 &37-75\\
39-84 &102-167 &47-95 &52-95  &103-159 &61-94& 55-98&104-161 &56-90\\
 63-118 &102-185 &44-110 & 45-95 &106-167 &44-108 &57-113&111-199 &46-83\\
62-116 &112-162 &63-109 &64-121 &130-180 &52-98 &67-122 &136-201 &62-95\\
55-97 &103-161 &56-84 &52-104 &90-177 &48-107 &59-101 &125-192 &54-92\\
58-109 &116-168 &26-109 &54-104 &97-182 &53-120 & 50-111&98-157 &61-108\\
57-101 &124-226 &49-88&47-108 &98-160 &54-78 &59-90 &120-180 &75-124\\
60-107 &97-154 &53-103 &54-104 &100-161 &58-99 &47-86 &87-150 &47-86\\
90-127 &159-214 &59-78 &77-158 &141-256 &70-132 &70-118 &138-221 &55-89\\
62-107 &108-147 &63-115 &50-95 &87-152 &55-80 &65-117 &115-196 &47-83\\
53-105 &120-188 &70-105 &42-86 &99-172 &56-103 &54-100 &95-166 &40-80\\
57-95 &113-176 &71-121 &45-107 &92-172 &56-97 &46-103 &114-186 &68-91\\
45-91 &83-140 &37-86 &100-136  &145-210 &62-100 & & &\\
\hline
\end{tabular}}
\vspace*{0.5cm}
\end{center}
\end{table}

\newpage
\begin{table}[!h]
\caption{\small Experimental results for the estimation of the linear
models}\label{table-sim}
\vspace{0.5cm}
\begin{center}
\begin{tabular}{ c l  c  c  c}
\hline
 \textbf{Model}  & $\widehat{\nu}$  $\setminus$ $n$ & $30$  & $100$ & $500$ \\
\hline
                 &              & $\widehat{E}(\widehat{\nu})$ \
$\widehat{MSE}(\widehat{\nu})$ & $\widehat{E}(\widehat{\nu})$ \
$\widehat{MSE}(\widehat{\nu})$ &  $\widehat{E}(\widehat{\nu})$ \
$\widehat{MSE}(\widehat{\nu})$ \\ 
\hline 
$M_1$ & $\widehat{b_1}$    & 1.9732 0.0042 &1.9858 0.0008  &1.9933 0.0001\\
      & $\widehat{b_2}$    &-4.9627  0.0056 &-4.9799 0.0013 &-4.9909 0.0002  \\
      & $\widehat{b_3}$    &-0.9809 0.0115 &-0.9926 0.0070 &-0.9961 0.0001  \\
 \hline
$M_2$ & $\widehat{b^1}$ &-2.0005  0.0097 &-1.9997  0.0026 &-2.0004 0.0005  \\
      & $\widehat{b^2}$ &1.9651 0.0052 &1.9809   0.0013 & 1.9911 0.0003 \\
      & $\widehat{b^3}$ &0.9302 0.0230 &0.9588   0.0060 &0.9816  0.0011 \\
      & $\widehat{b^4}$ &0.4991 0.0044 &0.5004  0.0012  &0.5000 0.0002\\
\hline
$M_3$  & $\widehat{b_1^1}$ &-2.0014 0.0114 &-2.0004   0.0026 &-2.0002 0.0005  \\
       & $\widehat{b_2^1}$ &5.0017 0.0465 &5.0007  0.0108 & 5.0001 0.0020 \\
       & $\widehat{b_3^1}$  &1.0002  0.0111 &1.0001  0.0027 &1.0000 0.0005  \\
       & $\widehat{b_1^2}$  &1.9738  0.0082 &1.9837   0.0019 &1.9920 0.0003 \\
       & $\widehat{b_2^2}$  &1.9763 0.0100 &1.9853 0.0020  &1.9920 0.0004\\
       & $\widehat{b_3^2}$  &0.9722 0.0082 &0.9841 0.0018 &0.9918 0.0003\\
       & $\widehat{b_1^3}$   &0.9576 0.0413 &0.9691 0.0090 &0.9855 0.0015\\
       & $\widehat{b_2^3}$   &0.9097 0.0737 &0.9429 0.0171 &0.9717 0.0030 \\
       & $\widehat{b_3^3}$   &2.9588 0.0410 &2.9709 0.0087 &2.9842 0.0015\\
       & $\widehat{b_1^4}$   &0.4996 0.0054 &0.5003  0.0013 &0.5001 0.0002 \\
       & $\widehat{b_2^4}$    &0.9992 0.0060 &1.0002 0.0014 &1.0002 0.0003\\
       & $\widehat{b_3^4}$   &-2.9994 0.0053 &-2.9995  0.0013 &-3.0002 0.0003\\
          \hline
\end{tabular}
\end{center}
\end{table}

\newpage
{\bf List of Figure Captions}

\begin{itemize}
\item Figure 1: $\Gamma_G$ for the sample data in Table 1
\item Figure 2: Box plot of the LS estimators for Model $M_1$, $n=30$ (left); $n=100$ (right)
\item Figure 3: Box plot of the LS estimators for Model $M_2$, $n=30$ (left); $n=100$ (right)
\item Figure 4: Box plot of the LS estimators for Model $M_3$, $n=30$ (left); $n=100$ (right)

\begin{algorithm}[h]
\caption{}
\begin{algorithmic}
\STATE  {\bf STEP 1:} Compute the global minimum of $g$, $\widehat{\nu}=S_2^{-1} z_2$, with 
$z_2=(\widehat{\sigma}_{\emph{\textbf{x}}^S,\emph{\textbf{y}}}    , 
\widehat{\sigma}_{\emph{\textbf{x}}^C,\emph{\textbf{y}}})^t$ and $S_2$ the
sample covariance matrix of  $(\emph{\textbf{x}}^S,\emph{\textbf{x}}^C)$.

\noindent \textbf{If} $\widehat{\nu} \in \Gamma_G$, \textbf{then} $\widehat{\nu}$ is the solution,
\textbf{else goto} Step 2. 

\STATE {\bf STEP 2:} Compute $r_0=\min_{i=1,\ldots,n} {\mathrm{spr}
\emph{\textbf{y}}_i}/{|\mathrm{mid} \emph{\textbf{x}}_i|}$ and identify the straight line $l_{(v)}$ in the set
$\{l_k: c= - u_k b + v_k\}_{k=1}^n$ such that $(0,r_0)\in l_{(v)}$. If there
exists more than one line in these conditions, then $l_{(v)}$ is the one for which the
value $- \mathrm{spr} \emph{\textbf{y}}_k / |\mathrm{mid} \emph{\textbf{x}}_k|$ is lowest.

\STATE  {\bf STEP 3:} Compute $s_0=\min_{i=1,\ldots,n} {\mathrm{spr}
\emph{\textbf{y}}_i}/{\mathrm{spr} \emph{\textbf{x}}_i}$ and identify the straight line $l_{(h)}$ in the set
$\{l_k: c= - u_k b + v_k\}_{k=1}^n$ such that $(s_0,0)\in l_{(h)}$. If there
exists more than one line in these conditions, then $l_{(h)}$ is the one for which the
value $- \mathrm{spr} \emph{\textbf{y}}_k / |\mathrm{mid} \emph{\textbf{x}}_k|$ is greatest.

\STATE  {ºbf STEP 4:} Let $R=\{l_{(v)}\}$, $C=\{0,s_0\}$, $D=\{(v),(h)\}$, $j=1$ and
$l_{(j)}=l_{(v)}$.

\textbf{If} $(v)=(h)$, \textbf{then} redefine $R=\{l^1\}, C=\{x^0,x^1\}$, let $t=1$ and \textbf{goto}
Step 8 \textbf{else goto} Step 5.

\STATE  {\bf STEP 5:} Compute $(b_{(j,h)},c_{(j,h)})$ the intersection point of the
lines $l_{(j)}$ and $l_{(h)}$. 

\noindent Check if $(b_{(j,h)},c_{(j,h)}) \in \mathrm{fr}(\Gamma_G)$, through
the conditions
\begin{itemize}
\item[i)] $b_{(j,h)} \in [0,s_0]$, and
\item[ii)]  $c_{(j,h)} = \min\{-u_k b_{(j,h)}+v_k : k=1,\ldots,n\}$.
\end{itemize}

\noindent \textbf{If} $(b_{(j,h)},c_{(j,h)}) \in \mathrm{fr}(\Gamma_G)$, \textbf{goto} Step 7
\textbf{else goto} Step 6. 

\STATE  {\bf STEP 6:} Compute $(b_{(j,k)},c_{(j,k)})$ the intersection points of
$l_{(j)}$ and each line in $\{l_k: c= - u_k b + v_k\}_{k=1}^n$ such that
$k\notin D$. Take the line $l_{k^*}$ such that $(b_{(j,k^*)},c_{(j,k^*)}) \in
\mathrm{fr}(\Gamma_G)$ (verifying the corresponding conditions i) and ii) shown
in Step 5). If there exists more than one line in these conditions, choose as
$l_{k^*}$ the one for which the value $- \mathrm{spr} \emph{\textbf{y}}_{k^*} / |\mathrm{mid} \emph{\textbf{x}}_{k^*}|$ is
lowest.

\noindent Let $R=R\cup \{l_{k*}\}$, $C=C\cup \{b_{(j,k*)}\}$, $D=D\cup \{k*\}$,
$j=j+1$, $l_{(j)}=l_{k*}$, and \textbf{goto} Step 5. 

\STATE {\bf STEP 7:} Let $R=R\cup \{l_{(h)}\}$ and $C=C\cup \{b_{(j,h)}\}$. 

\noindent Redefine
$R=\{l_{(v)},l_{k^{*}_{1}},l_{k^{*}_{2}},\ldots,l_{k^{*}_{p}},l_{(h)}\}$ as
$\{l^1,l^2,l^3,\ldots,l^{t-1},l^t\}$, and $C=\{0,b_{(1,k^{*}_{1})},
b_{(k^{*}_{1},k^{*}_{2})},\ldots, b_{(k^{*}_{p},h)},s_0\}$ as
$\{x^0,x^1,x^2,\ldots, x^{t-1},x^t\}$. \textbf{Goto} Step 8.
\end{algorithmic}
\end{algorithm}

\begin{algorithm}
\begin{algorithmic}
\STATE  {\bf STEP 8:} For $i=1,\ldots,t$, compute the local minimum of $g$ over the
segment corresponding to the line $l^i$ on $[x^{i-1},x^i]$, given by the
analytic expressions
$$\left\{\begin{array}{l} b_{*}^{i}=\max\big\{x^{i-1},\min\{b^{i},x^{i}\}\big\}
\\ 
c_{*}^{i}=- u_{i} b_{*}^{i}+v_{i}\end{array}\right.$$
\noindent where $\displaystyle b^{i}=\frac{u_i v_i \widehat{\sigma}_{\emph{\textbf{x}}^C}^2 -
v_i \widehat{\sigma}_{\emph{\textbf{x}}^S,\emph{\textbf{x}}^C}-u_i \widehat{\sigma}_{\emph{\textbf{x}}^C,\emph{\textbf{y}}} +
\widehat{\sigma}_{\emph{\textbf{x}}^S,\emph{\textbf{y}}}}{\widehat{\sigma}_{\emph{\textbf{x}}^S}^2 + u_i^2
\widehat{\sigma}_{\emph{\textbf{x}}^C}^2 -2 u_i \widehat{\sigma}_{\emph{\textbf{x}}^S,\emph{\textbf{x}}^C}}.$

\noindent Compute $g(b_{*}^{i},c_{*}^{i})$.

\noindent Take $(b_{L_2},c_{L_2})$ the point in
$\{(b_{*}^{i},c_{*}^{i})\}_{i=1}^{t}$ for which the value $g(b_{*}^{i},c_{*}^{i})$ is
lowest. Note that $(b_{L_2},c_{L_2})$ is the local minimum of $g$ over $L_2$. 

\STATE {\bf STEP 9:} Compute $(b_{L_1},c_{L_1})$ the local minimum of $g$ over $L_1$, given by the analytic expressions
$$\left\{\begin{array}{l} b_{L_1}=0 \\ 
c_{L_1}=\max\Big\{0,\min\bigg\{\displaystyle\frac{\widehat{\sigma}_{\emph{\textbf{x}}^C,\emph{\textbf{y}}}}{\widehat{\sigma}_{
\emph{\textbf{x}}^C}^2},r_0\bigg\}\Big\} \end{array}\right.$$
\noindent Compute $g(b_{L_1},c_{L_1})$.

\STATE {\bf STEP 10:} Compute $(b_{L_3},c_{L_3})$ the local minimum of $g$ over $L_3$, given by the analytic expressions
$$\left\{\begin{array}{l}
b_{L_3}=\max\Big\{0,\min\bigg\{\displaystyle\frac{\widehat{\sigma}_{\emph{\textbf{x}}^S,\emph{\textbf{y}}}}{\widehat{\sigma}_{
\emph{\textbf{x}}^S}^2},s_0\bigg\}\Big\}\\ 
c_{L_3}=0 \end{array}\right.$$
\noindent Compute $g(b_{L_3},c_{L_3})$.

\STATE {\bf STEP 11:} Take $(b^{*},c^{*})$ the point in
$\{(b_{L_j},c_{L_j})\}_{j=1}^{3}$ whose value  $g(b_{L_j},c_{L_j})$ is lowest.
Note that $(b^{*},c^{*})$ is the local minimum of $g$ on $\mathrm{fr}(S_G)$.

\end{algorithmic}
\end{algorithm}
\end{itemize}

\newpage
\end{document}